\documentclass{article}
\usepackage{amsmath}
\usepackage{amssymb}
\usepackage{mathtools}
\usepackage{enumerate}
\usepackage{graphicx}
\usepackage{float}
\usepackage{overpic}
\usepackage{relsize}
\usepackage{geometry}
\usepackage{algorithm}
\usepackage{algorithmic}
\usepackage{multirow}
\usepackage{stmaryrd}
\usepackage{tikz}

\newcommand\bC{\mathbb{C}}

\newcommand{\ap}[1]{\lfloor#1\rfloor}
\newcommand{\app}[1]{\llfloor#1\rrfloor}

\begin{document}
\title{Recursive Sweeping Preconditioner for the 3D Helmholtz Equation}
\author{Fei Liu$^\sharp$ and Lexing Ying$^{\dagger\sharp}$\\
  $\dagger$ Department of Mathematics, Stanford University\\
  $\sharp$ Institute for Computational and Mathematical Engineering, Stanford University
}
\date{Feb 2015}
\maketitle

\begin{abstract}
  This paper introduces the recursive sweeping preconditioner for the
  numerical solution of the Helmholtz equation in 3D. This is based on
  the earlier work of the sweeping preconditioner with the moving
  perfectly matched layers (PMLs). The key idea is to apply the
  sweeping preconditioner recursively to the quasi-2D auxiliary
  problems introduced in the 3D sweeping preconditioner. Compared to
  the non-recursive 3D sweeping preconditioner, the setup cost of this
  new approach drops from $O(N^{4/3})$ to $O(N)$, the application cost
  per iteration drops from $O(N\log N)$ to $O(N)$, and the iteration
  count only increases mildly when combined with the standard GMRES
  solver. Several numerical examples are tested and the results are
  compared with the non-recursive sweeping preconditioner to
  demonstrate the efficiency of the new approach.
\end{abstract}



\section{Introduction}
\label{sec:intro}
Let the domain of interest be the unit cube $D=(0,1)^3$ for
simplicity. The time-independent wave
field $u(x)$ satisfies the Helmholtz equation
\begin{equation}
  \Delta u(x)+\dfrac{\omega^2}{c^2(x)}u(x)=f(x),\quad \forall x\in D,
\end{equation}
where $\omega$ is the angular frequency, $c(x)$ is the velocity field
with a bound $c_{\min}\le c(x)\le c_{\max}$ where $c_{\min}$ and
$c_{\max}$ are assumed to be of $\Theta(1)$, and $f(x)$ is the
time-independent external force. The typical boundary conditions for
this problem are approximations of the Sommerfeld radiation condition,
which means that the wave is absorbed by the boundary and there is no
reflection coming from it. Other boundary conditions, such as the
Dirichlet boundary condition, can also be specified on part of the
boundary depending on the modeling setup.

In this setting, $\omega/(2\pi)$ is the typical wave number of the
problem and $\lambda=2\pi/\omega$ is the typical wavelength. For most
applications, the Helmholtz equation is discretized with at least a
few number of points (typically 4 to 20) per wavelength. So the number
of points $n$ in each direction is at least proportional to
$\omega$. As a result, the total degree of freedom
$N=n^3=\Omega(\omega^3)$ can be very large for high frequency 3D
problems. In addition, the corresponding discrete system is highly
indefinite and the standard iterative solvers and/or preconditioners
are no longer efficient for such problems. These together make the
problem challenging for numerical solution. We refer to the review
article \cite{why} by Ernst and Gander for more details on this.

Recently in \cite{sweeppml}, Engquist and Ying developed a sweeping
preconditioner using the moving perfectly matched layers (PMLs) and
obtained essentially linear solve times for 3D high frequency
Helmholtz equations. A key step of that approach is to approximate the
3D problem with a sequence of $O(n)$ PML-padded auxiliary quasi-2D
problems, each of which can be solved efficiently with sparse direct
method such as the nested dissection algorithm. As an extension, this
paper applies the sweeping idea recursively to further reduce each
auxiliary quasi-2D problem into a sequence of PML-padded quasi-1D
problems, each of which can be solved easily with the sparse LDU
factorization for banded systems. As a result, the setup cost of the
preconditioner improves from $O(N^{4/3})$ to $O(N)$ and the
application cost reduces from $O(N\log N)$ to $O(N)$.


There has been a vast literature on iterative methods and
preconditioners for the Helmholtz equation and we refer to the review
articles \cite{advances} by Erlangga and \cite{why} by Ernst and
Gander for a rather complete discussion. The discussion here only
touches on the methods that share similarity with the sweeping
preconditioners. The analytic ILU factorization (AILU) \cite{ailu} is
the first to use incomplete LDU factorizations for preconditioning the
Helmholtz equation. Compared to the moving PML sweeping
preconditioner, the method uses the absorbing boundary condition
(ABC), which is less effective compared to the PML, and hence the
iteration count grows much more rapidly.

Since the sweeping preconditioners \cite{sweephmf,sweeppml} were
proposed, there have been a number of exciting developments for the
numerical solutions of the high frequency Helmholtz equation,
including but not limited to
\cite{stolk2013domaindecomp,parallelsweep,sweepem,sweepemfem,sweepspectral,vion2014doublesweep,chen2013sourcetrans,chen2013sourcetrans2,demanet}. In
\cite{stolk2013domaindecomp}, Stolk proposed a domain decomposition
algorithm that utilizes suitable transmission conditions based on the
PMLs between the subdomains to achieve a near-linear cost. In
\cite{parallelsweep}, Poulson et al discussed a parallel version of
the moving PML sweeping preconditioner to deal with large scale
problems from applications such as seismic inversion. In
\cite{sweepem,sweepemfem,sweepspectral}, Tsuji and co-authors extended
the moving PML sweeping preconditioner method to other time-harmonic
wave equations and more general numerical discretization schemes. In
\cite{vion2014doublesweep}, Vion and Geuzaine proposed a double sweep
algorithm, studied several implementations of the absorbing boundary
conditions, and compared their numerical performance. Finally in
\cite{chen2013sourcetrans,chen2013sourcetrans2}, Chen and Xiang
introduced a sweeping-style domain decomposition method where the
emphasis was on the source transferring between the adjacent
subdomains. In \cite{demanet}, Zepeda-N{\'u}{\~n}ez and Demanet
developed a novel parallel domain decomposition method that uses
transmission conditions to define explicitly the up- and down-going
waves.

The rest of the paper is organized as follows. We first state the
problem and the discretization used in Section
\ref{sec:problem}. Section \ref{sec:review} reviews the non-recursive
moving PML sweeping preconditioner proposed in \cite{sweeppml}.
Section \ref{sec:recur} discusses in detail the recursive sweeping
preconditioner. Numerical results are presented in Section
\ref{sec:num}. Finally, the conclusion and some future directions are
provided in Section \ref{sec:conclu}.

\section{Problem Formulation}
\label{sec:problem}

Following \cite{sweeppml}, we assume that the perfectly matched layer
(PML) \cite{berenger1994pml,chew1994pml,johnson2008pmlnotes} is
utilized at part of the boundary where the Sommerfeld radiation
condition is specified. The sweeping preconditioner in \cite{sweeppml}
requires that at least one of the six faces of the domain
$D=(0,1)^3$ is specified with the PML boundary condition. As we shall
see soon, the recursive sweeping preconditioner instead requires the PML
condition to be specified at least at two non-parallel faces. Without
loss of generality, we assume that it is specified at $x_2=0$ and
$x_3=0$. There is no restriction on the type of boundary conditions
specified on the other four faces. However, to simplify the
discussion, we assume that the Dirichlet condition is used. The PML
boundary condition introduces auxiliary functions
\begin{gather*}
  \sigma(x)=
  \begin{dcases}
    \dfrac{C}{\eta}\left(\dfrac{x-\eta}{\eta}\right)^2,&x\in[0,\eta],\\
    0,& x\in(\eta,1],
  \end{dcases}
\end{gather*}
and
\begin{gather*}
  s(x)=\left(1+i\dfrac{\sigma(x)}{\omega}\right)^{-1},
  \quad
  s_1(x)\equiv 1,
  \quad
  s_2(x)=s(x_2),
  \quad
  s_3(x)=s(x_3),
\end{gather*}
where $C$ is an appropriate positive constant independent of $\omega$,
and $\eta$ is the PML width, which is typically around one
wavelength. The Helmholtz equation with PML is
\begin{equation}
  \begin{dcases}
    \left((s_1\partial_1)(s_1\partial_1)+(s_2\partial_2)(s_2\partial_2)+(s_3\partial_3)(s_3\partial_3)+\dfrac{\omega^2}{c^2(x)}\right)u(x)=f(x), & \forall x\in D=(0,1)^3,\\
    u(x)=0, & \forall x\in \partial D.
  \end{dcases}
  \label{eq:helmpml}
\end{equation}
It is typically assumed that that the support of $f(x)$ is in
$(0,1)\times(\eta,1)\times(\eta,1)$, which means that the force is not
located in the PML region.  The cube $[0,1]^3$ is discretized with a
Cartesian grid where the grid size is $h=\frac{1}{n+1}$ and $n$ is
proportional to $\omega$. The set of all the interior points of the
grid is given by
\begin{gather*}
  P=\{p_{i,j,k}=(ih,jh,kh):1\le i,j,k \le n\},
\end{gather*}
and the degree of freedom is $N=n^3$.

\begin{figure}[h!]
  \begin{tikzpicture}[scale=0.7]
    \fill[gray](-6cm,-6cm)rectangle(0cm,0cm);
    \fill[gray,xscale=0.5,yslant=0.5](0,-6cm)rectangle(6cm,0cm);
    \fill[gray,yscale=0.5,xslant=0.5](-6cm,0cm)rectangle(0cm,0.5cm)rectangle(-0.5cm,6cm);
    \draw[step=0.25cm](-6cm,-6cm)grid(0cm,0cm)(-3cm,-3cm)node{PML}[->](0cm,-6cm)--(-6.5cm,-6cm)node[below]{$x_3$};
    \draw[->](0cm,-6cm)--(0cm,0.5cm)node[above]{$x_1$};
    \draw[step=0.25cm,xscale=0.5,yslant=0.5](0,-6cm)grid(6cm,0cm)(3cm,-3cm)node[xscale=0.5,yslant=0.5]{PML}[->](0,-6cm)--(6.5cm,-6cm)node[right,xscale=0.5,yslant=0.5]{$x_2$};
    \draw[step=0.25cm,yscale=0.5,xslant=0.5](-6cm,0cm)grid(0cm,6cm)(-3cm,3cm)node[yscale=0.5,xslant=0.5]{zero Dirichlet};
  \end{tikzpicture}
  \begin{tikzpicture}
    \fill[gray](0cm,6cm)rectangle(0.5cm,0cm)rectangle(6cm,0.5cm);
    \draw[step=0.25cm](0cm,0cm)grid(6cm,6cm);
    \draw[->](0,0)--(6.5cm,0)node[below]{$x_2$};
    \draw[->](0,0)--(0,6.5cm)node[left]{$x_3$};
    \draw(3cm,0cm)node[below]{PML}(0cm,3cm)node[above,rotate=90]{PML}(3cm,6cm)node[above]{zero Dirichlet}(6cm,3cm)node[rotate=-90,above]{zero Dirichlet};
  \end{tikzpicture}
  \caption{The domain of interest. Left is a 3D view of the
    domain. Right is an $x_2$-$x_3$ cross section view, where each
    cell stands for a 1D column. The gray area stands for the PML
    region. }
  \label{fig:1}
\end{figure}
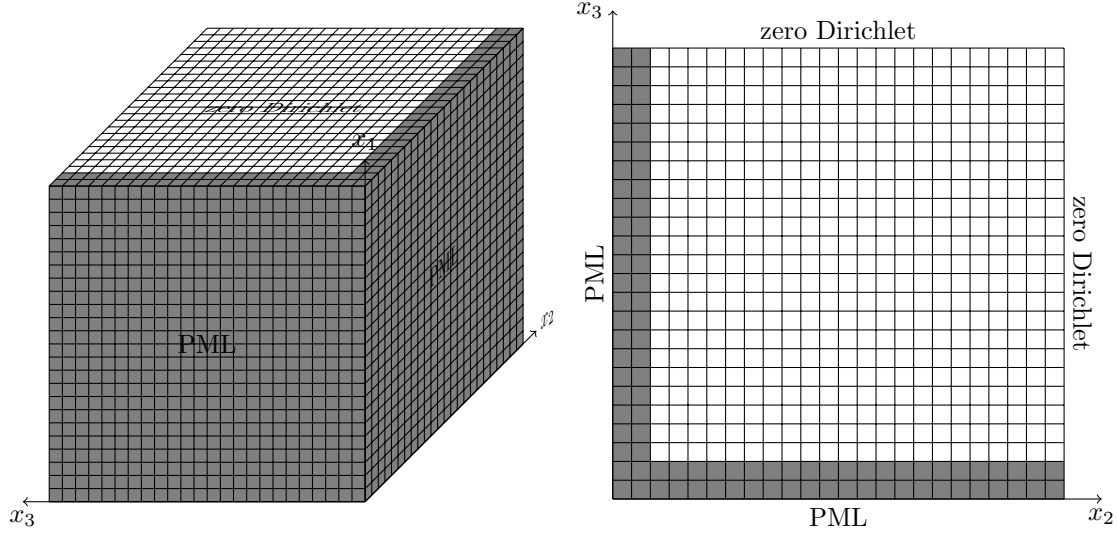

Applying the standard 7-point finite difference stencil results in the
discretized system
\begin{equation}
  \label{eq:system}
  \begin{gathered}
    \dfrac{(s_1)_{i,j,k}}{h}\left(\dfrac{(s_1)_{i+1/2,j,k}}{h}(u_{i+1,j,k}-u_{i,j,k})-\dfrac{(s_1)_{i-1/2,j,k}}{h}(u_{i,j,k}-u_{i-1,j,k})\right)\\ +\dfrac{(s_2)_{i,j,k}}{h}\left(\dfrac{(s_2)_{i,j+1/2,k}}{h}(u_{i,j+1,k}-u_{i,j,k})-\dfrac{(s_2)_{i,j-1/2,k}}{h}(u_{i,j,k}-u_{i,j-1,k})\right)\\ +\dfrac{(s_3)_{i,j,k}}{h}\left(\dfrac{(s_3)_{i,j,k+1/2}}{h}(u_{i,j,k+1}-u_{i,j,k})-\dfrac{(s_3)_{i,j,k-1/2}}{h}(u_{i,j,k}-u_{i,j,k-1})\right)\\ +\left(\dfrac{\omega^2}{c^2}\right)_{i,j,k}u_{i,j,k}=f_{i,j,k},
    \quad \forall 1\le i,j,k\le n,
  \end{gathered}
\end{equation}
where the subscript $(i,j,k)$ means that the corresponding function is
evaluated at the point $p_{i,j,k}=(ih,jh,kh)$ and the definition of the points here extends to half integers as well. The computational task
is to solve \eqref{eq:system} efficiently. We note that, unlike the
symmetric version adopted in \cite{sweephmf,sweeppml}, here the
nonsymmetric version of the equation is used. Figure \ref{fig:1}
provides an illustration of the computational domain and the
discretization grid.

\section{Review of the Sweeping Preconditioner with Moving PML}
\label{sec:review}
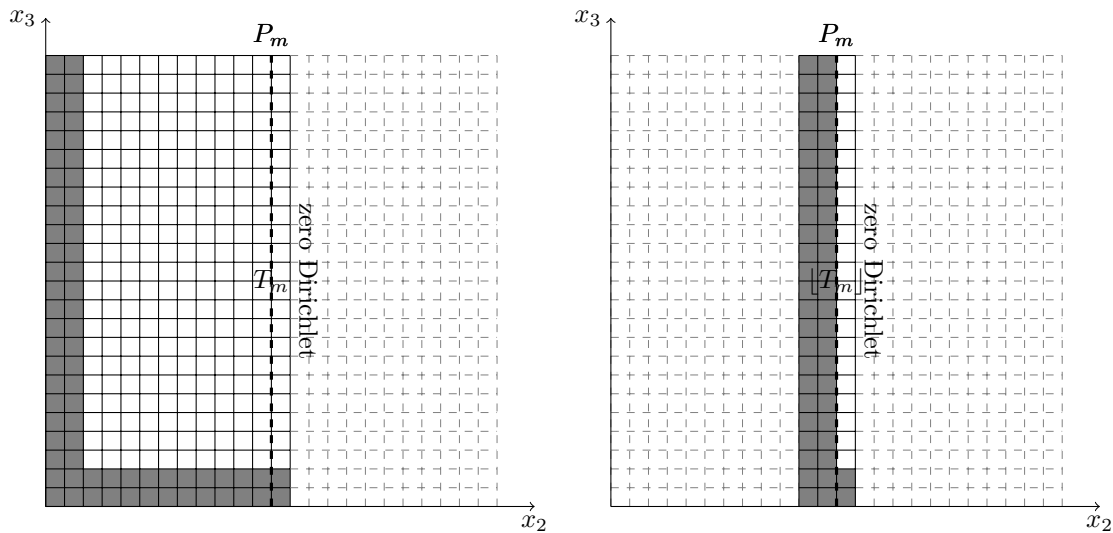
\begin{figure}[h!]
  \begin{tikzpicture}
    \fill[gray](0,0)rectangle(0.5cm,6cm);
    \fill[gray](0,0)rectangle(3.25cm,0.5cm);
    \draw[step=0.25cm](0,0)grid(3.25cm,6cm);
    \draw[step=0.25cm,dashed,gray,very thin](3.25cm,0)grid(6cm,6cm);
    \draw[->](0,0)--(6.5cm,0)node[below]{$x_2$};
    \draw[->](0,0)--(0,6.5cm)node[left]{$x_3$};
    \draw[dashed,very thick](3cm,0)--(3cm,6cm);
    \draw(3cm,6cm)node[anchor=south]{$P_{m}$}(3.25cm,3cm)node[above,rotate=-90]{zero Dirichlet};
    \draw(3cm,6cm)node[anchor=south]{$P_{m}$}(3cm,3cm)node{$T_m$};
  \end{tikzpicture}
  \begin{tikzpicture}
    \fill[gray](2.5cm,0)rectangle(3cm,6cm);
    \fill[gray](2.5cm,0)rectangle(3.25cm,0.5cm);
    \draw[step=0.25cm,dashed,gray,very thin](0,0)grid(2.5cm,6cm);
    \draw[step=0.25cm,dashed,gray,very thin](3.25cm,0)grid(6cm,6cm);
    \draw[step=0.25cm](2.5cm,0)grid(3.25cm,6cm);
    \draw[->](0,0)--(6.5cm,0)node[below]{$x_2$};
    \draw[->](0,0)--(0,6.5cm)node[left]{$x_3$};
    \draw[dashed,very thick](3cm,0)--(3cm,6cm);
    \draw(3cm,6cm)node[anchor=south]{$P_{m}$}(3.25cm,3cm)node[above,rotate=-90]{zero Dirichlet};
    \draw(3cm,6cm)node[anchor=south]{$P_{m}$}(3cm,3cm)node{$\ap{T_m}$};
    \draw(2.5cm,0cm)--(2.5cm,6cm);
  \end{tikzpicture}
  \caption{Left: $T_m$ is the restriction to $P_m$ (the dashed
    grid) of the half space Green's function on the solid
    grid. Right: By moving the PML right next to the layer $P_m$,
    the operator $T_m$ is approximated by solving the equation on a
    much smaller grid. }
  \label{fig:approx}
\end{figure}

This section gives a brief review of the non-recursive moving PML
sweeping preconditioner proposed in \cite{sweeppml} for
completeness. More details can be found in the original paper
\cite{sweeppml}. The starting point of the sweeping preconditioner is
a block LDU factorization called the sweeping factorization. To build
this factorization, the algorithm sweeps along the $x_3$ direction
starting from the face $x_3=0$. The unknowns with subscript index
$(i,j,k)$ are ordered with column-major order, i.e., first dimension
$1$, then dimension $2$, and finally dimension $3$. We
define the vectors
\begin{gather*}
  u=[u_{1,1,1},\dots,u_{n,1,1},\dots,u_{n,n,1},\dots,u_{n,n,n}]^T,\\
  f=[f_{1,1,1},\dots,f_{n,1,1},\dots,f_{n,n,1},\dots,f_{n,n,n}]^T.
\end{gather*}
By introducing
\[
  P_m=\{p_{1,1,m},\dots,p_{n,1,m},\dots,p_{n,n,m}\}
\]
as the points on the $m$-th plane and also
\begin{gather*}
  u_{:,:,m}=[u_{1,1,m},\dots,u_{n,1,m},\dots,u_{n,n,m}]^T,\\
  f_{:,:,m}=[f_{1,1,m},\dots,f_{n,1,m},\dots,f_{n,n,m}]^T,
\end{gather*}
one can write the system \eqref{eq:system} compactly as $Au=f$ with
the following block form
\begin{equation}
  \begin{bmatrix}
    A_{1,1}&A_{1,2}\\
    A_{2,1}&A_{2,2}&\ddots\\
    &\ddots&\ddots&A_{n-1,n}\\
    &&A_{n,n-1}&A_{n,n}
  \end{bmatrix}
  \begin{bmatrix}
    u_{:,:,1}\\u_{:,:,2}\\\vdots\\u_{:,:,n}
  \end{bmatrix}
  =
  \begin{bmatrix}
    f_{:,:,1}\\f_{:,:,2}\\\vdots\\f_{:,:,n}
  \end{bmatrix}.
  \label{eq:Aufblock}
\end{equation}
By defining $S_k$ and $T_k$ recursively via
\begin{gather*}
  S_1=A_{1,1},\quad T_1=S_1^{-1},\\
  S_m=A_{m,m}-A_{m,m-1}T_{m-1}A_{m-1,m},\quad T_m=S_m^{-1},\quad m=2,\dots,n,
\end{gather*}
the standard block LDU factorization of the block tridiagonal matrix
$A$ is
\begin{gather*}
  A=L_1\dots L_{n-1}
  \begin{bmatrix}
    S_1\\
    &\ddots\\
    &&S_n
  \end{bmatrix}
  U_{n-1}\dots U_1,
\end{gather*}
where $L_m$ and $U_m$ are the corresponding unit lower and upper
triangular matrices with the only non-zero off-diagonal blocks
\begin{gather*}
  L_m(P_{m+1},P_m)=A_{m+1,m}T_m, \quad U_m(P_{m},P_{m+1})=T_mA_{m,m+1}, \quad m=1,\dots, n-1.
\end{gather*}
It is not difficult to see that computing this factorization takes
$O(N^{7/3})$ steps. Once it is available, $u$ can be computed in
$O(N^{5/3})$ steps by
\begin{gather*}
  u=
  \begin{bmatrix}
    u_{:,:,1}\\\vdots\\u_{:,:,n}
  \end{bmatrix}
  =A^{-1}f=U_1^{-1}\dots U_{n-1}^{-1}
  \begin{bmatrix}
    T_1\\
    &\ddots\\
    &&T_n
  \end{bmatrix}
  L_{n-1}^{-1}\dots L_1^{-1}f
\end{gather*}
The main disadvantage of the above algorithm is, $S_m$ and $T_m$ are
in general dense matrices of size $n^2\times n^2$ so the
corresponding dense linear algebra operations are expensive. The
sweeping preconditioner overcomes this difficulty by approximating
$T_m$ efficiently for $P_m$ with $mh\in(\eta,1]$, i.e., for $P_m$ not in
the PML region at the face $x_3=0$. The key point is to consider the physical meaning of
$T_m$. From now on let us assume $\eta=bh$ which implies that there
are $b$ layers in the PML region at $x_3=0$. Restricting the
factorization to the upper-left $m\times m$ block of $A$ where $m=b+1,\dots,n$ gives
\begin{gather*}
  \begin{bmatrix}
    A_{1,1}&A_{1,2}\\
    A_{2,1}&A_{2,2}&\ddots\\
    &\ddots&\ddots&A_{m-1,m}\\
    &&A_{m,m-1}&A_{m,m}
      \end{bmatrix}=L_1\dots L_{m-1}
  \begin{bmatrix}
    S_1\\&S_2\\
    &&\ddots\\
    &&&S_m
  \end{bmatrix}
  U_{m-1}\dots U_1,
\end{gather*}
where $L_t$ and $U_t$ are redefined by restricting to their upper left
$m\times m$ blocks. Inverting both sides leads to
\begin{gather*}
  \begin{bmatrix}
    A_{1,1}&A_{1,2}\\
    A_{2,1}&A_{2,2}&\ddots\\
    &\ddots&\ddots&A_{m-1,m}\\
    &&A_{m,m-1}&A_{m,m}
  \end{bmatrix}^{-1}
  =
  U_1^{-1}\dots U_{m-1}^{-1}
  \begin{bmatrix}
    T_1\\&T_2\\
    &&\ddots\\
    &&&T_m
  \end{bmatrix}
  L_{m-1}^{-1}\dots L_1^{-1}.
\end{gather*}
The left-hand side is the discrete half-space Green's function with
Dirichlet zero boundary condition at $x_3=(m+1)h$ and a
straightforward calculation shows that the lower-right block of the
right-hand side is $T_m$. Therefore, $T_m$ is the discrete half-space
Green's function restricted to the $m$-th layer. Note that, the PML at
$x_3=0$ is used to simulate an absorbing boundary condition. If we
assume that there is little reflection during the transmission of the
wave, we can approximate $T_m$ by placing the PML right next to the
$m$-th layer since the domain of interest is only the $m$-th layer (see
Figure \ref{fig:approx}). This is the key idea of the moving PML
sweeping preconditioner, where the operator $T_m$ is numerically
approximated by putting the PML right next to the domain of interest
and solving a much smaller system to save the computational cost.

More precisely, we introduce an auxiliary problem on the domain
$D_m=[0,1] \times [0,1] \times[(m-b)h,(m+1)h]: $
\begin{gather*}
  \begin{dcases}
    \left((s_1\partial_1)(s_1\partial_1)+(s_2\partial_2)(s_2\partial_2)+(s_3^m\partial_3)(s_3^m\partial_3)+\dfrac{\omega^2}{c^2(x)}\right)v(x)=g(x),& \forall x\in D_m,\\
    v(x)=0, & \forall x\in \partial D_m,
  \end{dcases}
\end{gather*}
where $s_3^m(x)=s(x_3-(m-b)h)$. The domain $D_m$ is discretized with
the partial grid
\begin{gather*}
  P_{(m-b+1):m}:=\{P_t:m-b+1\le t\le m\}.
\end{gather*}
Applying the same central finite difference scheme gives rise to the
corresponding discretized system, denoted as
\begin{gather*}
  H_m v=g, \quad m=b+1,\ldots,n.
\end{gather*}
To approximate $T_m$, we numerically define operator
$\ap{T_m}:\alpha\in\bC^{n^2}\to \beta\in\bC^{n^2}$ by the
following procedure:
\begin{enumerate}
\item
  Introduce a vector $g$ defined on $P_{(m-b+1):m}$ by setting
  $\alpha$ to the layer $P_m$ and zero everywhere else.
\item
  Solve the discretized auxiliary problem $H_m v=g$ on $P_{(m-b+1):m}$
  with $g$ from step 1.
\item
  Set $\beta$ as the restriction on $P_m$ of the solution $v$ from
  step 2.
\end{enumerate}
The discretized system is a quasi-2D system as $b$ is
typically a small constant, so the system can be solved efficiently by
the nested dissection method
\cite{george1973nested,duff1983multifrontal,liu1992multifrontal}.

The first $b$ layers, which are in the PML region of the original
problem \eqref{eq:helmpml}, need to be handled with a slight
difference. Define
\begin{gather*}
  u_{:,:,1:b}=[u_{:,:,1}^T,\dots,u_{:,:,b}^T]^T,\\
  f_{:,:,1:b}=[f_{:,:,1}^T,\dots,f_{:,:,b}^T]^T.
\end{gather*}
Then the system $Au=f$ can be written as
\begin{gather*}
  \begin{bmatrix}
    A_{1:b,1:b}&A_{1:b,b+1}\\
    A_{b+1,1:b}&A_{b+1,b+1}&\ddots\\
    &\ddots&\ddots&A_{n-1,n}\\
    &&A_{n,n-1}&A_{n,n}
  \end{bmatrix}
  \begin{bmatrix}
    u_{:,:,1:b}\\u_{:,:,b+1}\\\vdots\\u_{:,:,n}
  \end{bmatrix}
  =
  \begin{bmatrix}
    f_{:,:,1:b}\\f_{:,:,b+1}\\\vdots\\f_{:,:,n}
  \end{bmatrix}.
\end{gather*}
For the first $b$ layers, we simply define $\ap{T_{1:b}}$ as the
inverse operator of $H_b:=A_{1:b,1:b}$. However, it is essential that $\ap{T_{1:b}}$ is
stored in a factorized form by applying the nested dissection method
to $H_b$, since $H_bv=g$ is also a quasi-2D problem.

Based on the above discussion, the setup algorithm of the moving PML
sweeping preconditioner is given in Algorithm \ref{alg:sweepsetup}.
\begin{algorithm}[h!]
  \caption{Construction of the moving PML sweeping preconditioner of the system \eqref{eq:system}. Complexity
    $=O(b^3n^4)=O(b^3N^{4/3})$. }
  \label{alg:sweepsetup}
  \begin{algorithmic}
    \STATE
    Construct the nested dissection factorization of $H_b$, which defines $\ap{T_{1:b}}$.
    \FOR {$m=b+1,\dots,n$}
    \STATE Construct the nested dissection factorization of $H_m$, which defines $\ap{T_m}$.
    \ENDFOR
  \end{algorithmic}
\end{algorithm}

Once the factorization is completed, $\ap{T_{1:b}}$ and $\ap{T_m}$ can be applied using
the nested dissection factorization. The application process of the sweeping
preconditioner is given in Algorithm \ref{alg:sweepsolve}.
\begin{algorithm}[h!]
  \caption{Computation of $u\approx A^{-1} f$ using the factorization from Algorithm \ref{alg:sweepsetup}. Complexity $=O(b^2n^3\log n)=O(b^2 N\log N)$. }
  \label{alg:sweepsolve}
  \begin{algorithmic}
    \STATE
    $u_{:,:,1:b}=\ap{T_{1:b}} f_{:,:,1:b}$\\
    $u_{:,:,b+1}=\ap{T_{b+1}}(f_{:,:,b+1}-A_{b+1,1:b}u_{:,:,1:b})$
    \FOR {$m=b+1,\dots,n-1$}
    \STATE
    $u_{:,:,m+1}=\ap{T_{m+1}}(f_{:,:,m+1}-A_{m+1,m}u_{:,:,m})$
    \ENDFOR
    \FOR{$m=n-1,\dots,b+1$}
    \STATE
    $u_{:,:,m}=u_{:,:,m}-\ap{T_m}(A_{m,m+1}u_{:,:,m+1})$
    \ENDFOR\\
    $u_{:,:,1:b}=u_{:,:,1:b}-\ap{T_{1:b}}(A_{1:b,b+1}u_{:,:,b+1})$
  \end{algorithmic}
\end{algorithm}

\section{Recursive Sweeping Preconditioner}
\label{sec:recur}

\newcommand{\td}{\widetilde}

Recall that the PML is also applied to the face $x_2=0$. Therefore,
each quasi-2D auxiliary problem is itself a discretization of the
Helmholtz equation with the PML specified on one side. Following the
treatment in \cite{sweeppml} for the 2D Helmholtz equation, it is natural
to apply the same sweeping idea once again along the $x_2$ direction,
instead of the nested dissection algorithm used in the previous
section.

\subsection{Inner sweeping}

\begin{figure}[h!]
  \begin{tikzpicture}
  \fill[gray](2.5,0)rectangle(3cm,3.25cm);
  \fill[gray](3cm,0)rectangle(3.25cm,0.5cm);
  \draw[step=0.25cm](2.5cm,0)grid(3.25cm,3.25cm);
  \draw(2.5cm,0cm)--(2.5cm,3.25cm);
  \draw[step=0.25cm,dashed,gray,very thin](0,0)grid(6cm,6cm);
  \draw[->](0,0)--(6.5cm,0)node[below]{$x_2$};
  \draw[->](0,0)--(0,6.5cm)node[left]{$x_3$};
  \draw[dashed,very thick](2.5cm,3cm)--(3.25cm,3cm);
  \draw(2.5cm,3cm)node[left]{$\td{P}_m$};
  \draw(3cm,3cm)node{$\td{T}_m$};
  \end{tikzpicture}
  \begin{tikzpicture}
  \fill[gray](2.5,2.5cm)rectangle(3cm,3.25cm);
  \fill[gray](3cm,2.5cm)rectangle(3.25cm,3cm);
  \draw[step=0.25cm](2.5cm,2.5cm)grid(3.25cm,3.25cm);
  \draw(2.5cm,3.25cm)--(2.5cm,2.5cm)--(3.25cm,2.5cm);
  \draw[step=0.25cm,dashed,gray,very thin](0,0)grid(6cm,6cm);
  \draw[->](0,0)--(6.5cm,0)node[below]{$x_2$};
  \draw[->](0,0)--(0,6.5cm)node[left]{$x_3$};
  \draw[dashed,very thick](2.5cm,3cm)--(3.25cm,3cm);
  \draw(2.5cm,3cm)node[left]{$\td{P}_m$};
  \draw(3cm,3cm)node{$\ap{\td{T}_m}$};
  \end{tikzpicture}
  \caption{Left: $\td{T}_m$ is the restriction to $\td{P}_m$ (the dashed grid) of the
    Green's function on the quasi-2D solid grid. Right: By moving the
    PML right next to $\td{P}_m$, the operator $\td{T}_m$ is
    approximated by solving the problem on a quasi-1D grid. }
  \label{fig:recurapprox}
\end{figure}
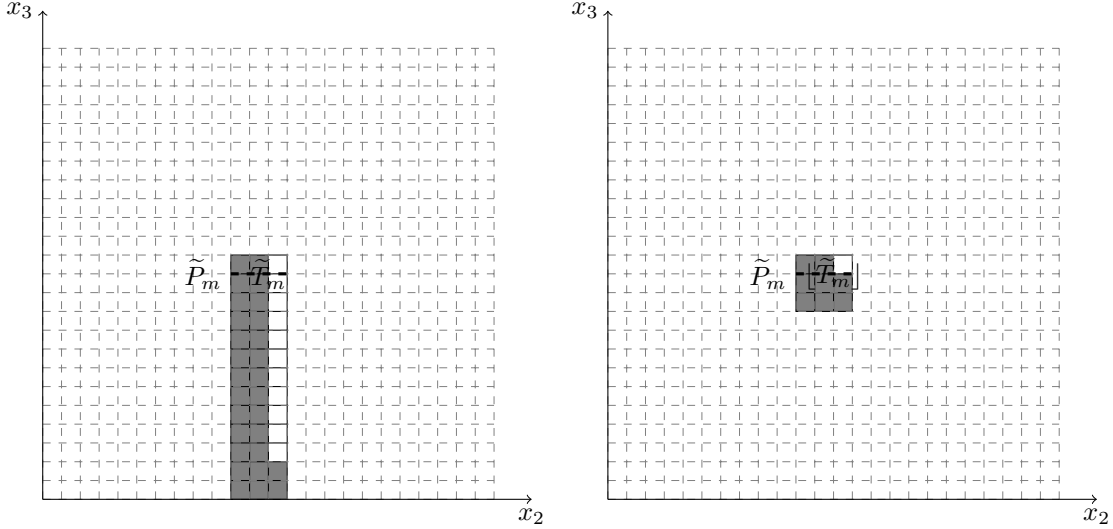

Recall that the quasi-2D subproblems of the non-recursive sweeping
preconditioners are $H_m v=g,m=b,\dots,n$. Since they have essentially
the same structure, it is sufficient to consider a single system
$\td{A} v = g$ where $\td{A}$ can be anyone of $H_m$.
Here the accent mark is to emphasize that the problem under
consideration is quasi-2D. To formalize the sweeping preconditioner
along the $x_2$ direction, we define, up to a translation,
\[
\td{P}=\{p_{i,j,k}=(ih,jh,kh):1\le i,j \le n, 1\le k\le b\},
\]
to be the discretization grid. For each $m=1,\ldots,n$, let
\begin{gather*}
  \td{P}_m=\{p_{1,m,1},\dots,p_{1,m,b},\dots,p_{n,m,b}\}, \\
  v_{:,m,:}=[v_{1,m,1},\dots,v_{1,m,b},\dots,v_{n,m,b}]^T,\\
  g_{:,m,:}=[g_{1,m,1},\dots,g_{1,m,b},\dots,g_{n,m,b}]^T.
\end{gather*}
For the first $b$ layers in the $x_2$ direction, we also define
\begin{gather*}
  \td{P}_{1:b}=\{\td{P}_1,\dots,\td{P}_b\},\\
  v_{:,1:b,:}=[v_{:,1,:}^T,\dots,v_{:,b,:}^T]^T, \\
  g_{:,1:b,:}=[g_{:,1,:}^T,\dots,g_{:,b,:}^T]^T.
\end{gather*}
In this section, we reorder the vectors $v,g$ by grouping the 3rd
dimension first and applying the column-major ordering to dimensions 1
and 2:
\begin{gather*}
  v=[v_{:,1,:}^T,\dots,v_{:,n,:}^T]^T, \\
  g=[g_{:,1,:}^T,\dots,g_{:,n,:}^T]^T.
\end{gather*}
With this ordering, the corresponding system $\td{A} v=g$ is written
as
\begin{gather*}
  \begin{bmatrix}
    \td{A}_{1:b,1:b} & \td{A}_{1:b,b+1}\\
    \td{A}_{b+1,1:b} & \td{A}_{b+1,b+1} & \ddots\\
    &\ddots&\ddots & \td{A}_{n-1,n}\\
    && \td{A}_{n,n-1} & \td{A}_{n,n}
  \end{bmatrix}
  \begin{bmatrix}
    v_{:,1:b,:}\\v_{:,b+1,:}\\\vdots\\v_{:,n,:}
  \end{bmatrix}
  =
  \begin{bmatrix}
    g_{:,1:b,:}\\g_{:,b+1,:}\\\vdots\\g_{:,n,:}
  \end{bmatrix}.
\end{gather*}
For the block LDU factorization of $\td{A}$, we define
\begin{gather*}
  \td{S}_{1:b}=\td{A}_{{1:b},{1:b}}, \quad \td{T}_{1:b}=\td{S}_{1:b}^{-1}, \\
  \td{S}_{b+1}=\td{A}_{b+1,b+1}-\td{A}_{b+1,{1:b}}\td{T}_{1:b}\td{A}_{{1:b},b+1},\quad \td{T}_{b+1}=\td{S}_{b+1}^{-1},\\
  \td{S}_m=\td{A}_{m,m}-\td{A}_{m,m-1}\td{T}_{m-1}\td{A}_{m-1,m},\quad \td{T}_m=\td{S}_m^{-1},  \quad m=b+2,\dots,n,
\end{gather*}
then $\td{A}$ can be factorized as
\begin{gather*}
  \td{A}=\td{L}_{1:b}\td{L}_{b+1}\dots \td{L}_{n-1}
  \begin{bmatrix}
    \td{S}_{1:b}\\
    &\td{S}_{b+1}\\
    &&\ddots\\
    &&&\td{S}_n
  \end{bmatrix}
  \td{U}_{n-1}\dots \td{U}_{b+1}\td{U}_{1:b},
\end{gather*}
where the non-zero off-diagonal blocks of the unit lower and upper triangular matrices $\td{L}_m$ and $\td{U}_m$ are given by
\begin{gather*}
  \td{L}_{1:b}(\td{P}_{b+1},\td{P}_{1:b})=\td{A}_{b+1,{1:b}}\td{T}_{1:b},
  \quad \td{U}_{1:b}(\td{P}_{1:b},\td{P}_{b+1})=\td{T}_{1:b}\td{A}_{{1:b},b+1},\\
  \td{L}_m(\td{P}_{m+1},\td{P}_m)=\td{A}_{m+1,m}\td{T}_m,
  \quad \td{U}_m(\td{P}_m,\td{P}_{m+1})=\td{T}_m\td{A}_{m,m+1}, \quad m=b+1,\dots,n-1.
\end{gather*}
Then the solution $v$ can be computed by
\begin{gather*}
  v=
  \begin{bmatrix}
    v_{:,{1:b},:}\\v_{:,b+1,:}\\\vdots\\v_{:,n,:}
  \end{bmatrix}
  =\td{A}^{-1}g=
  \td{U}_{1:b}^{-1}\td{U}_{b+1}^{-1}\dots \td{U}_{n-1}^{-1}
  \begin{bmatrix}
    \td{T}_{1:b}\\
    &\td{T}_{b+1}\\
    &&\ddots\\
    &&&\td{T}_n
  \end{bmatrix}
  \td{L}_{n-1}^{-1}\dots \td{L}_{b+1}^{-1}\td{L}_{1:b}^{-1}g.
\end{gather*}
By comparing the factorization of the upper-left $(m-b+1)\times(m-b+1)$
block of $\td{A}$, where $m=b+1,\dots,n$, we have
\begin{gather*}
  \begin{bmatrix}
    \td{A}_{{1:b},{1:b}}&\td{A}_{{1:b},b+1}\\
    \td{A}_{b+1,{1:b}}&\td{A}_{b+1,b+1}&\ddots\\
    &\ddots&\ddots&\td{A}_{m-1,m}\\
    &&\td{A}_{m,m-1}&\td{A}_{m,m}
  \end{bmatrix}
  =\td{L}_{1:b}\td{L}_{b+1}\dots \td{L}_{m-1}
  \begin{bmatrix}
    \td{S}_{1:b}\\
    &\td{S}_{b+1}\\
    &&\ddots\\
    &&&\td{S}_m
  \end{bmatrix}
  \td{U}_{m-1}\dots \td{U}_{b+1}\td{U}_{1:b},
\end{gather*}
where $\td{L}_t$ and $\td{U}_t$ are redefined as their
restrictions to their top-left $(m-b+1)\times(m-b+1)$
blocks. Inverting both sides gives
\begin{gather*}
  \begin{bmatrix}
    \td{A}_{{1:b},{1:b}}&\td{A}_{{1:b},b+1}\\
    \td{A}_{b+1,{1:b}}&\td{A}_{b+1,b+1}&\ddots\\
    &\ddots&\ddots&\td{A}_{m-1,m}\\
    &&\td{A}_{m,m-1}&\td{A}_{m,m}
  \end{bmatrix}^{-1}
  =
  \td{U}_{1:b}^{-1}\td{U}_{b+1}^{-1}\dots \td{U}_{m-1}^{-1}
  \begin{bmatrix}
    \td{T}_{1:b}\\
    &\td{T}_{b+1}\\
    &&\ddots\\
    &&&\td{T}_m
  \end{bmatrix}
  \td{L}_{m-1}^{-1}\dots \td{L}_{b+1}^{-1}\td{L}_{1:b}^{-1}.
\end{gather*}
Thus, by repeating the argument in Section \ref{sec:review}, the
matrix $\td{T}_m$ is the restriction to the layer $\td{P}_m$ of the
discrete half-space Green's function. It can be approximated by
$\ap{\td{T}_m}$, which is defined by solving a quasi-1D
problem obtained by placing a moving PML right next to $x_2=mh$ (see
Figure \ref{fig:recurapprox}). Each auxiliary quasi-1D problem in
this inner sweeping step can be solved by the sparse block LDU
factorization efficiently, with ordering the system by grouping dimension $3$ and $2$ first and dimension $1$ last.

More specifically, for each $m$, we introduce the auxiliary
problem on the domain
$\td{D}_m=[0,1]\times[(m-b)h,(m+1)h]\times[0,(b+1)h]$:
\begin{gather*}
  \begin{dcases}
    \left((s_1\partial_1)(s_1\partial_1)+(s_2^m\partial_2)(s_2^m\partial_2)+(s_3\partial_3)(s_3\partial_3)+\dfrac{\omega^2}{c^2(x)}\right)w(x)=q(x),&\forall x\in \td{D}_m,\\
    w(x)=0,&\forall x\in \partial \td{D}_m,
  \end{dcases}
\end{gather*}
where $s_2^m(x)=s(x_2-(m-b)h)$. The domain $\td{D}_m$ is discretized
with the grid
\begin{gather*}
  \td{P}_{(m-b+1):m}:=\{\td{P}_t:m-b+1\le t\le m\},
\end{gather*}
and the same central difference numerical scheme is used here. We
denote the corresponding discretized system as $\td{H}_m w=q$.
Similar to the process described in Section \ref{sec:review}, we
define the operator $\ap{\td{T}_m}:\alpha\in\bC^{nb}\to
\beta\in\bC^{nb}$ by the following procedure:
\begin{enumerate}
\item
  Introduce a vector $q$ defined on the grid $\td{P}_{(m-b+1):m}$ by
  setting $\alpha$ to the layer $\td{P}_m$ and zero everywhere else.
\item
  Solve the auxiliary quasi-1D problem $\td{H}_m w=q$ on
  $\td{P}_{(m-b+1):m}$ with $q$ from step 1.
\item
  Set $\beta$ as the restriction on $\td{P}_m$ of the solution $w$
  from step 2.
\end{enumerate}
For the first $b$ layers, $\ap{\td{T}_{1:b}}$ is simply defined as the
inverse operator of $\td{H}_b:=\td{A}_{1:b,1:b}$, which is essentially the same as $\td{T}_{1:b}$, but implemented by using the sparse
block LDU factorization of $\td{H}_b$. Summarizing all this, the setup
and application algorithm of the inner moving PML sweeping preconditioner are given in
Algorithms \ref{alg:recursetup} and \ref{alg:recursolve},
respectively.

\begin{algorithm}[h!]
  \caption{Construction of the inner moving PML sweeping preconditioner of the quasi-2D problem $\td{A}v=g$. Complexity $=O(b^6n^2)$. }
    \label{alg:recursetup}
  \begin{algorithmic}
    \STATE
    Construct the sparse block LDU factorization of $\td{H}_b$, which defines $\ap{\td{T}_{1:b}}$.
    \FOR {$m=b+1,\dots,n$}
    \STATE Construct the sparse block LDU factorization of $\td{H}_m$, which defines $\ap{\td{T}_m}$.
    \ENDFOR
  \end{algorithmic}
\end{algorithm}

\begin{algorithm}[h!]
  \caption{Computation of $v\approx \td{A}^{-1} g$ using the factorization
    from Algorithm \ref{alg:recursetup}. Complexity $=O(b^4n^2)$. }
  \label{alg:recursolve}
  \begin{algorithmic}
    \STATE
    $v_{:,{1:b},:}=\ap{\td{T}_{1:b}} g_{:,{1:b},:}$\\
    $v_{:,b+1,:}=\ap{\td{T}_{b+1}}(g_{:,b+1,:}-\td{A}_{b+1,{1:b}}v_{:,{1:b},:})$
    \FOR {$m=b+1,\dots,n-1$}
    \STATE
    $v_{:,m+1,:}=\ap{\td{T}_{m+1}}(g_{:,m+1,:}-\td{A}_{m+1,m}v_{:,m,:})$
    \ENDFOR
    \FOR{$m=n-1,\dots,b+1$}
    \STATE
    $v_{:,m,:}=v_{:,m,:}-\ap{\td{T}_m}(\td{A}_{m,m+1}v_{:,m+1,:})$
    \ENDFOR\\
    $v_{:,{1:b},:}=v_{:,{1:b},:}-\ap{\td{T}_{1:b}}(\td{A}_{{1:b},b+1}v_{:,b+1,:})$
  \end{algorithmic}
\end{algorithm}

\subsection{Putting together}

As we pointed out earlier, the matrix $\td{A}$ can be anyone of $H_m,m=b,\dots,n$, where Algorithms \ref{alg:recursetup} and
\ref{alg:recursolve} can be applied. Notice that solving the subproblems exactly with the nested dissection algorithm
results in the approximation $\ap{T_m}$ to $T_m$. This extra-level of
approximation defines a further approximation, which shall be denoted
by $\app{T_m}:\alpha\in\bC^{n^2}\to \beta\in\bC^{n^2}$ (to be precise, for the first $b$ layers, it is $\app{T_{1:b}}:\alpha\in\bC^{n^2b}\to \beta\in\bC^{n^2b}$). The steps for
carrying out $\app{T_m}$ are similar to the ones for $\ap{T_m}$ except
that one uses Algorithms \ref{alg:recursetup} and \ref{alg:recursolve}
to solve the quasi-2D problems approximately (instead of the nested
dissection method that solves them exactly).

Given all these preparations, the setup algorithm of the recursive
sweeping preconditioner can be summarized compactly in Algorithm
\ref{alg:allsetup} and the application algorithm is given in Algorithm
\ref{alg:allsolve}.
\begin{algorithm}[h!]
  \caption{Construction of the recursive moving PML sweeping preconditioner of
    the linear system \eqref{eq:system}. Complexity
    $=O(b^6n^3)=O(b^6N)$. }
  \label{alg:allsetup}
  \begin{algorithmic}
    \STATE Construct the inner moving PML sweeping preconditioner of $H_b$ by Algorithm
    \ref{alg:recursetup}. This gives $\app{T_{1:b}}$.
    \FOR {$m=b+1,\dots,n$}
    \STATE Construct the inner moving PML sweeping preconditioner of $H_m$ by Algorithm \ref{alg:recursetup}. This gives $\app{T_m}$.
    \ENDFOR
  \end{algorithmic}
\end{algorithm}
\begin{algorithm}[h!]
  \caption{Computation of $u\approx A^{-1} f$ using the factorization from Algorithm \ref{alg:allsetup}. Complexity $=O(b^4n^3)=O(b^4N)$. }
  \label{alg:allsolve}
  \begin{algorithmic}
    \STATE
    $u_{:,:,{1:b}}=\app{T_{1:b}}f_{:,:,{1:b}}$\\
    $u_{:,:,b+1}=\app{T_{b+1}}(f_{:,:,b+1}-A_{b+1,{1:b}}u_{:,:,{1:b}})$
    \FOR {$m=b+1,\dots,n-1$}
    \STATE
    $u_{:,:,m+1}=\app{T_{m+1}}(f_{:,:,m+1}-A_{m+1,m}u_{:,:,m})$
    \ENDFOR
    \FOR{$m=n-1,\dots,b+1$}
    \STATE
    $u_{:,:,m}=u_{:,:,m}-\app{T_m}(A_{m,m+1}u_{:,:,m+1})$
    \ENDFOR\\
    $u_{:,:,{1:b}}=u_{:,:,{1:b}}-\app{T_{1:b}}(A_{{1:b},b+1}u_{:,:,b+1})$
  \end{algorithmic}
\end{algorithm}

In the outer loop of Algorithm \ref{alg:allsolve}, the unknowns are
eliminated layer by layer in the $x_3$ direction. In the application of $\app{T_m}$, there is the inner loop in which the unknowns in each quasi-2D problem
are eliminated in the $x_2$ direction. The whole algorithm serves as a
preconditioner for the original linear system
\eqref{eq:system}. Notice that, in the recursive sweeping
preconditioner, the quasi-2D problems are solved only
approximately. Therefore, the overall accuracy might not be as good as
the non-recursive method. But as we will show in the next section, the
performance of the preconditioner is only mildly affected.

The above algorithms are described in a way to present the main ideas
clearly. In the actual implementations, a couple of modifications are
taken in order to maximize the efficiency:
\begin{enumerate}
\item
  For each auxiliary problem, both in the inner loop and the
  outer loop, several layers are processed together instead of one
  layer.
\item
  For the PML introduced in the auxiliary problems, the number of
  layers in the auxiliary PML region does not have to match the number
  of layers $b$ used for the boundary PML at $x_2=0$ and $x_3=0$. In
  fact, the thickness of the auxiliary PML is typically thinner for the
  sake of efficiency.
\item
  The problem we described above has zero Dirichlet boundary
  conditions on the other four faces of the cube. If instead, the PMLs
  are put on all the faces, then the sweeping preconditioner sweeps
  with two fronts from two opposite faces respectively and they meet in the middle
  with a subproblem with PML on both sides instead of only one side,
  as described in \cite{sweeppml}.
\end{enumerate}

\section{Numerical Results}
\label{sec:num}

In this section we test several numerical examples to illustrate the performance of the recursive sweeping preconditioner. All algorithms
are implemented in MATLAB and the tests are performed on a 2.0 GHz
computer with 256 GB memory. We use the GMRES algorithm as the
iterative solver with relative residual $10^{-3}$ and restart
value $40$ for the entire 3D system. The quasi-2D problems are solved approximately by applying the inner sweeping preconditioner only once for the sake of efficiency. The velocity fields and forces tested are kept the same with \cite{sweeppml} so that the results can be
compared easily. The PMLs are put on all six sides of the cube $[0,1]^3$
to simulate the Sommerfeld radiation condition.

\begin{figure}[h!]
  \centering
  \begin{overpic}
    [width=0.32\textwidth]{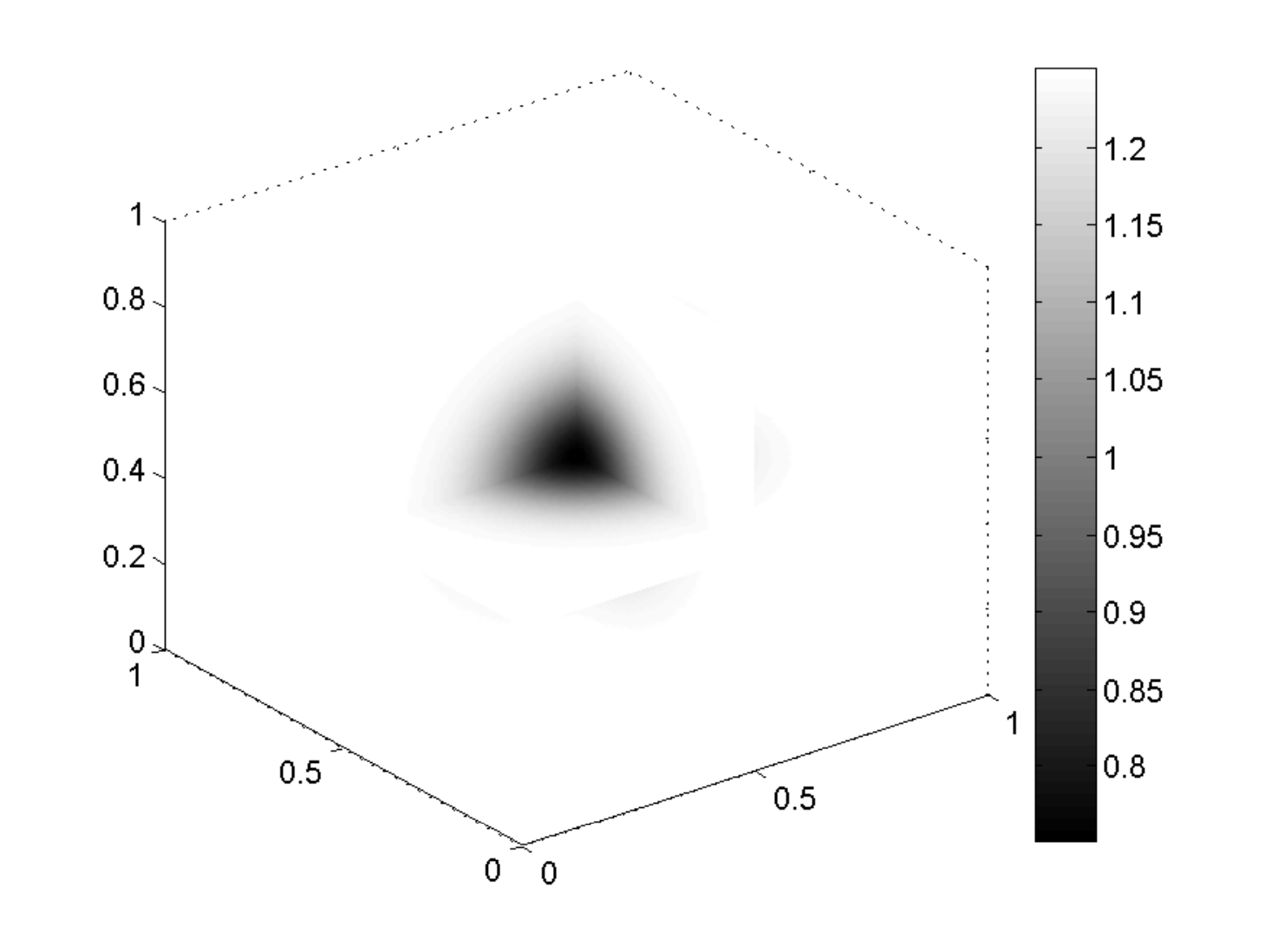}
    \put(41,-3){$(a)$}
  \end{overpic}
  \begin{overpic}
    [width=0.32\textwidth]{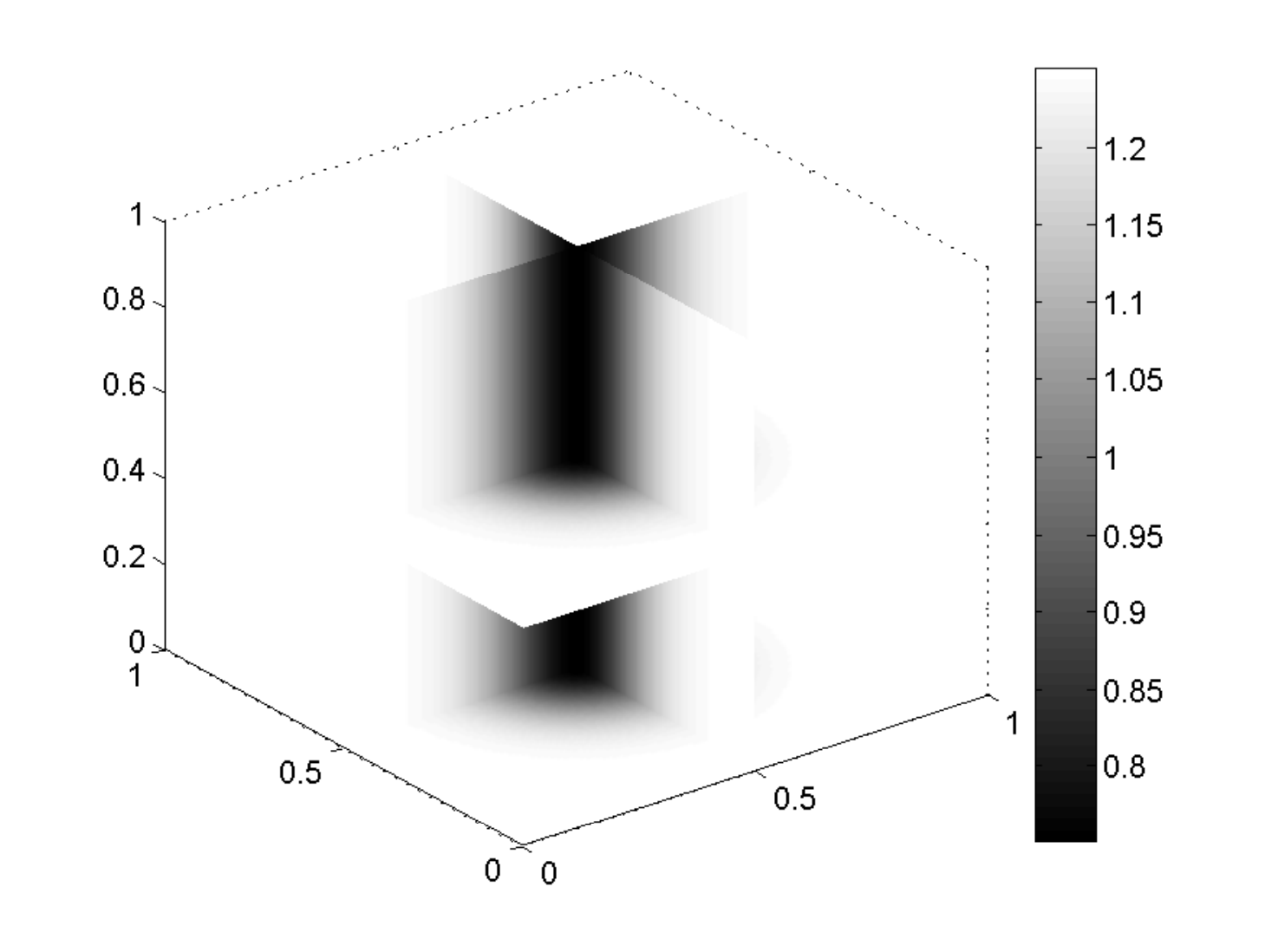}
    \put(41,-3){$(b)$}
  \end{overpic}
  \begin{overpic}
    [width=0.32\textwidth]{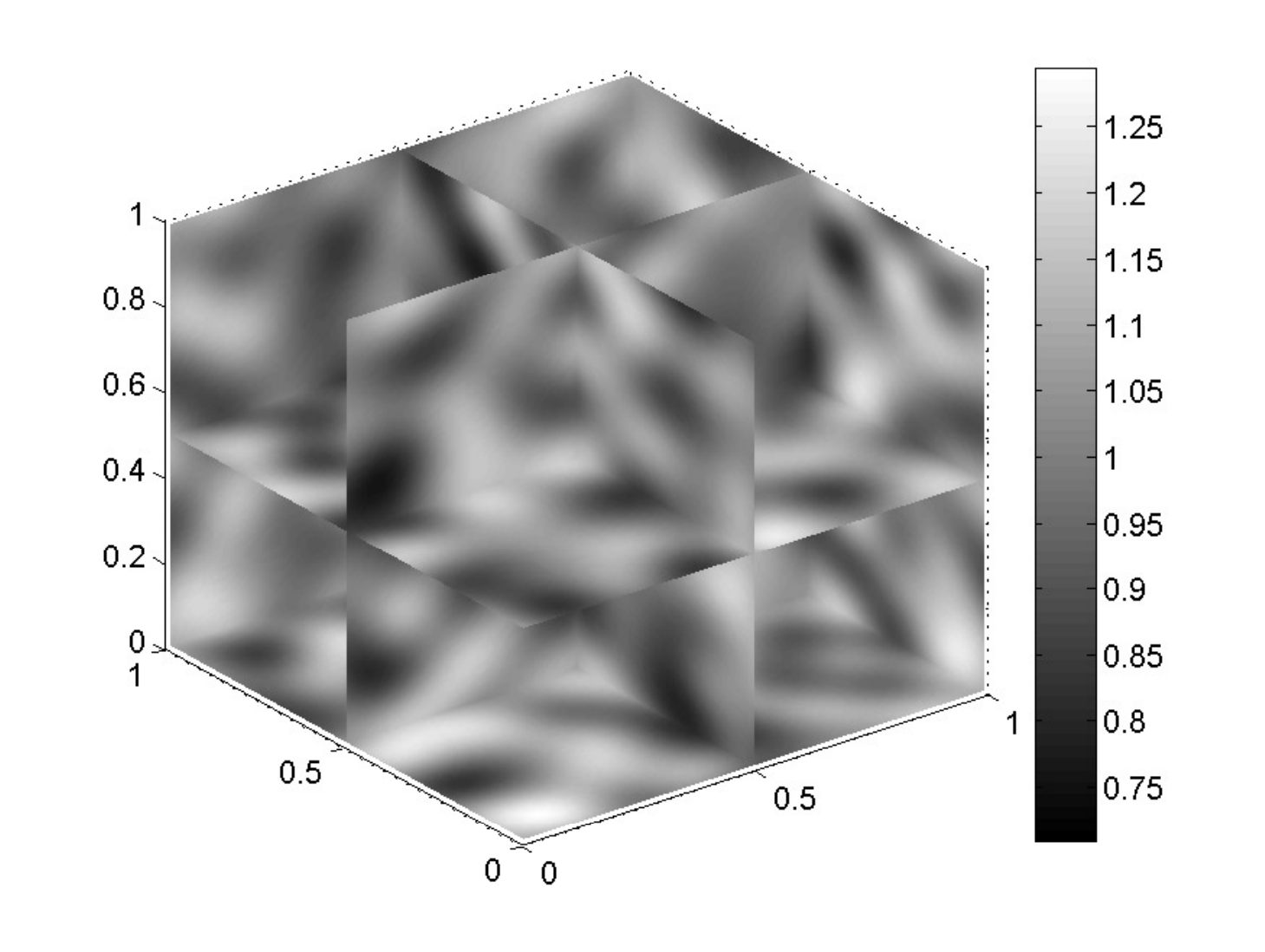}
    \put(41,-3){$(c)$}
  \end{overpic}
  \caption{The three velocity fields tested. }
  \label{fig:c}
\end{figure}

We test three velocity
fields (see Figure \ref{fig:c}):
\begin{enumerate}[(a)]
\item
  A converging lens with a Gaussian profile at the center of the domain.
\item
  A vertical waveguide with a Gaussian cross-section.
\item
  A random velocity field.
\end{enumerate}

For each velocity field, the tests are performed for two external
forces:
\begin{enumerate}[(a)]
\item
  A Gaussian point source centered at $(1/2,1/2,1/4)$.
\item
  A Gaussian wave packet with wavelength comparable to the typical
  wavelength of the domain. The packet centers at $(1/2,1/4,1/4)$ and
  points to the direction $(0,1/\sqrt{2},1/\sqrt{2})$.
\end{enumerate}
We vary the typical wave number $\omega/(2\pi)$, test the behavior of
the recursive preconditioner, and compare the results with the
non-recursive preconditioner.

In these tests, each wavelength is discretized with $q=8$ points. The
width of the PML at the boundary of the cube is $9h$, and the width of
the auxiliary PML for the middle layers is $5h$. The number of layers
processed in each auxiliary problem is $4$. The algorithm sweeps with
two fronts from $x_3=0$ and $x_3=1$ in the outer loop, and with two
fronts from $x_2=0$ and $x_2=1$ in the inner loop.

The results are reported in the following tables. $T_{\text{setup}}$
is the time used to construct the preconditioner in
seconds. $T_{\text{solve}}$ is the time used to solve the system in
the preconditioned GMRES solver in seconds and $N_{\text{iter}}$ is
the corresponding iteration number. ``NR'' stands for the original
non-recursive method while ``R'' stands for the recursive method
introduced in this paper. The ``ratio'' is the time cost of the recursive
method over the non-recursive method. The numerical implementation of
the non-recursive method is slightly improved as compared to
\cite{sweeppml}, by incorporating a more accurate PML
discretization. Therefore, the results here for the non-recursive
method are better compared to the ones in \cite{sweeppml}.

\begin{table}[h!]
  \centering
  \begin{overpic}
    [width=0.45\textwidth]{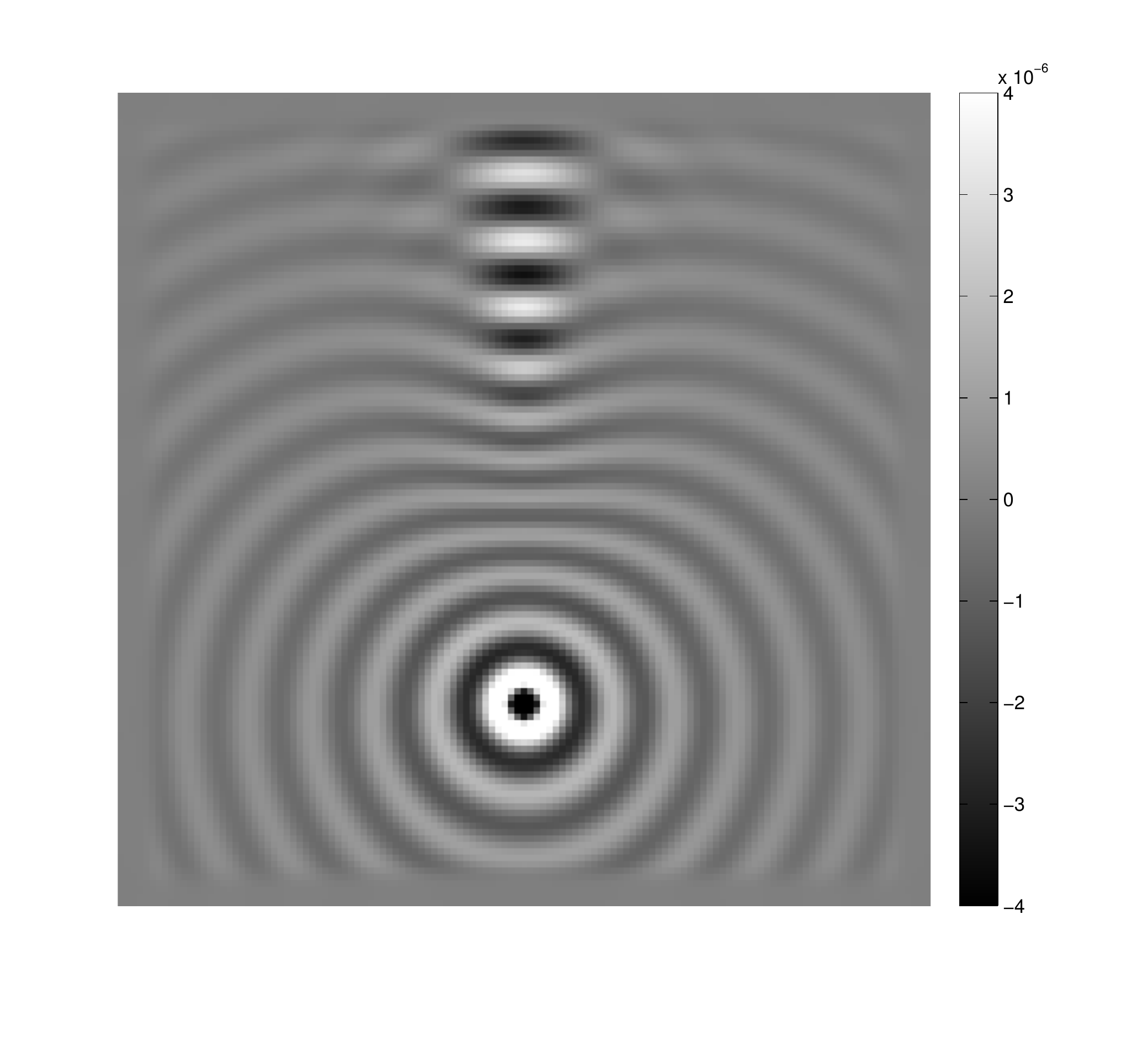}
    \put(42,6){$(a)$}
  \end{overpic}
  \begin{overpic}
    [width=0.45\textwidth]{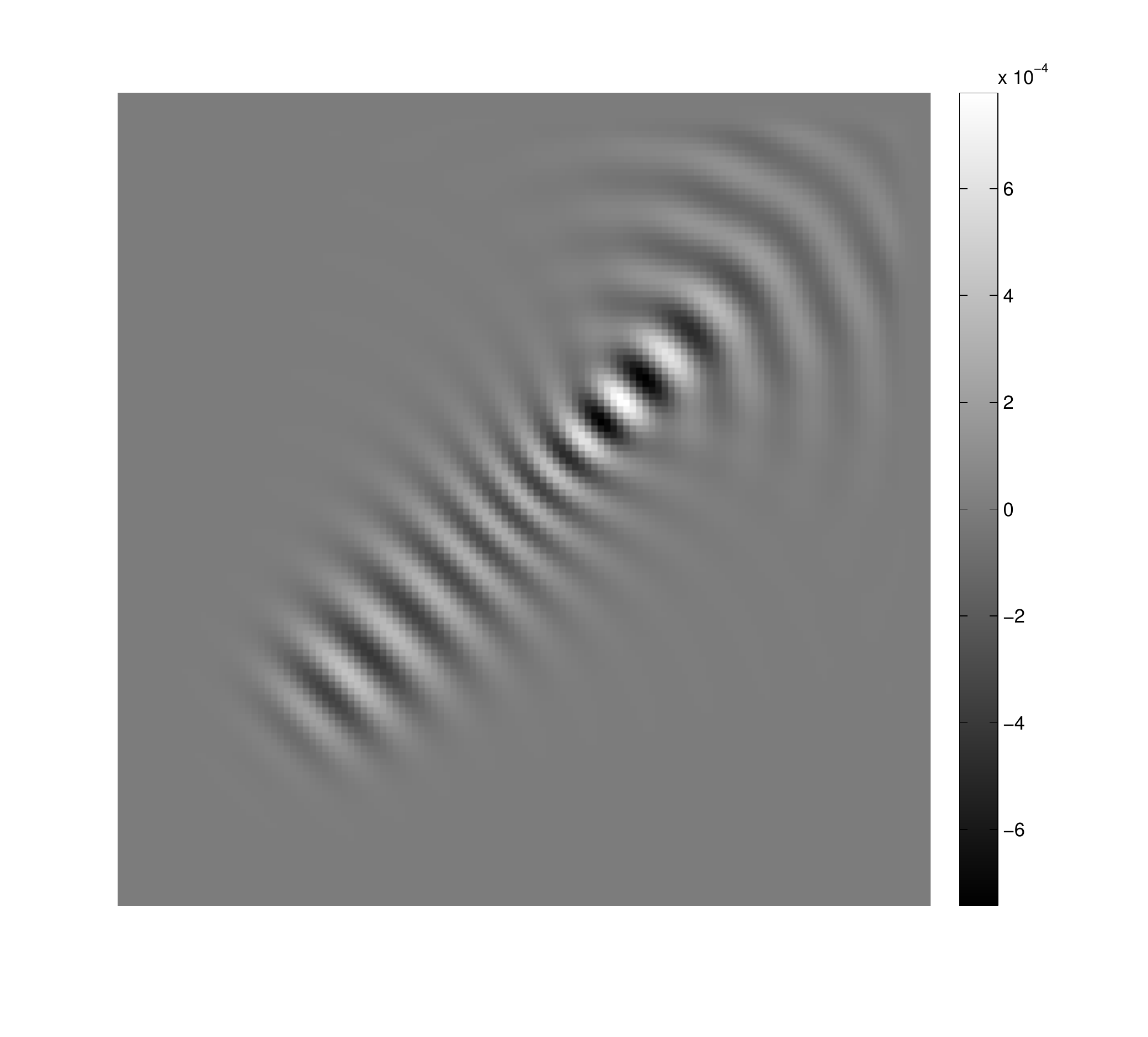}
    \put(42,6){$(b)$}
  \end{overpic}
  \begin{tabular}{ccc|ccc|c|cc|ccc}
    \hline
    &&&\multicolumn{3}{|c|}{$T_{\text{setup}}$}&&\multicolumn{2}{|c}{$N_{\text{iter}}$}&\multicolumn{3}{|c}{$T_{\text{solve}}$}\\
    \hline
    $\omega/(2\pi)$&$q$&$N$& NR & R &ratio&$f(x)$& NR & R & NR & R &ratio\\
     \hline
     \multirow{2}{*}{8}&\multirow{2}{*}{8}&\multirow{2}{*}{$63^3$}&\multirow{2}{*}{46.923}&\multirow{2}{*}{19.823}&\multirow{2}{*}{42\%}&$(a)$&3&4&12.313&16.355&133\%\\
     &&&&&&$(b)$&4&4&14.973&16.862&113\%\\
     \hline
     \multirow{2}{*}{16}&\multirow{2}{*}{8}&\multirow{2}{*}{$127^3$}&\multirow{2}{*}{537.12}&\multirow{2}{*}{180.99}&\multirow{2}{*}{34\%}&$(a)$&3&4&116.44&169.34&145\%\\
     &&&&&&$(b)$&4&4&150.67&168.65&112\%\\
     \hline
     \multirow{2}{*}{32}&\multirow{2}{*}{8}&\multirow{2}{*}{$255^3$}&\multirow{2}{*}{5927.0}&\multirow{2}{*}{1308.0}&\multirow{2}{*}{22\%}&$(a)$&4&5&1273.0&2039.8&160\%\\
     &&&&&&$(b)$&4&5&1312.1&2070.4&158\%\\
     \hline
  \end{tabular}
  \caption{Results for velocity field (a) in Figure
    \ref{fig:c}. Solutions with $\omega/(2\pi)=16$ at $x_1=0.5$ are
    presented. }
\end{table}

\begin{table}[h!]
  \centering
    \begin{overpic}
      [width=0.45\textwidth]{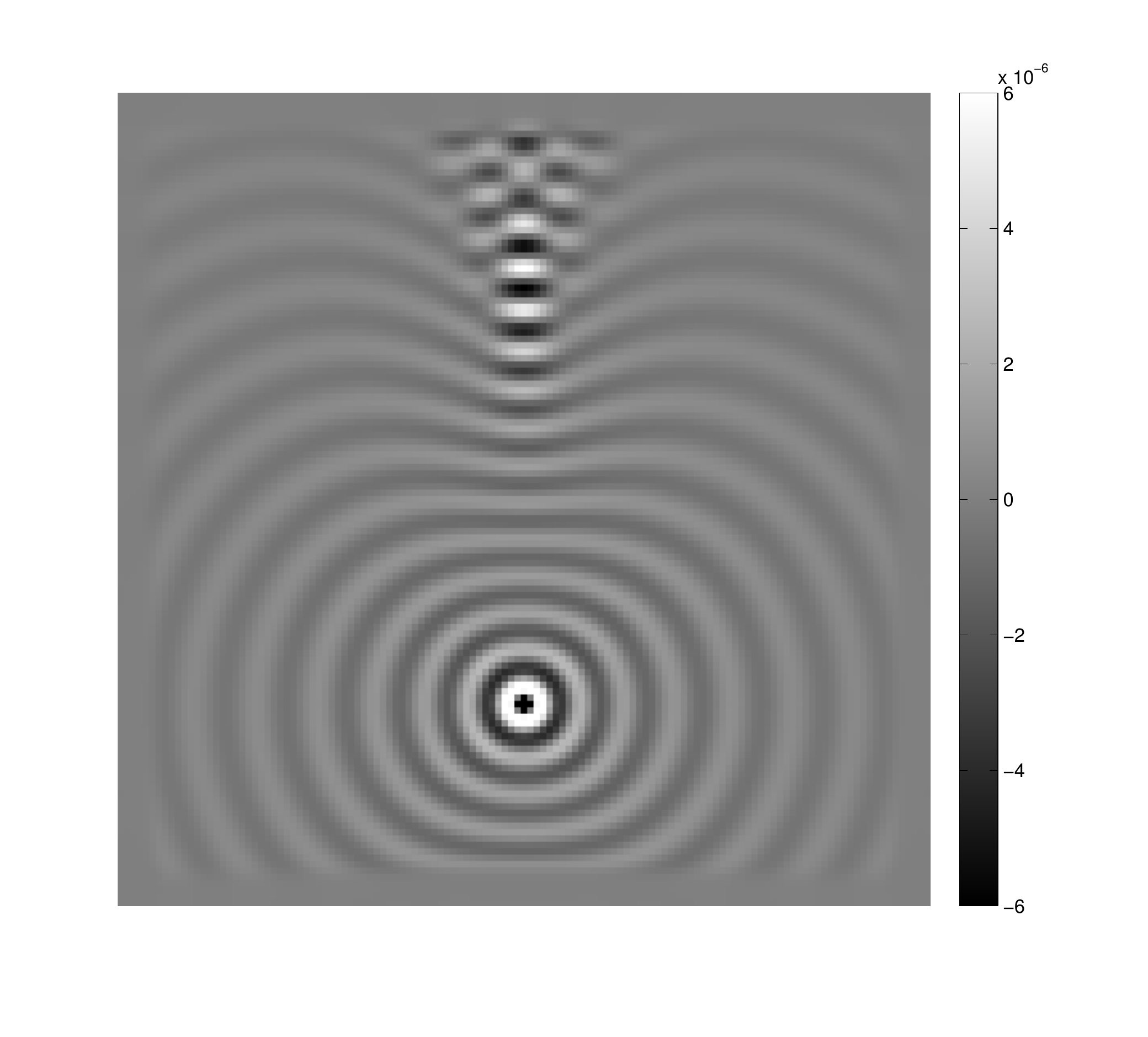}
      \put(42,6){$(a)$}
    \end{overpic}
    \begin{overpic}
      [width=0.45\textwidth]{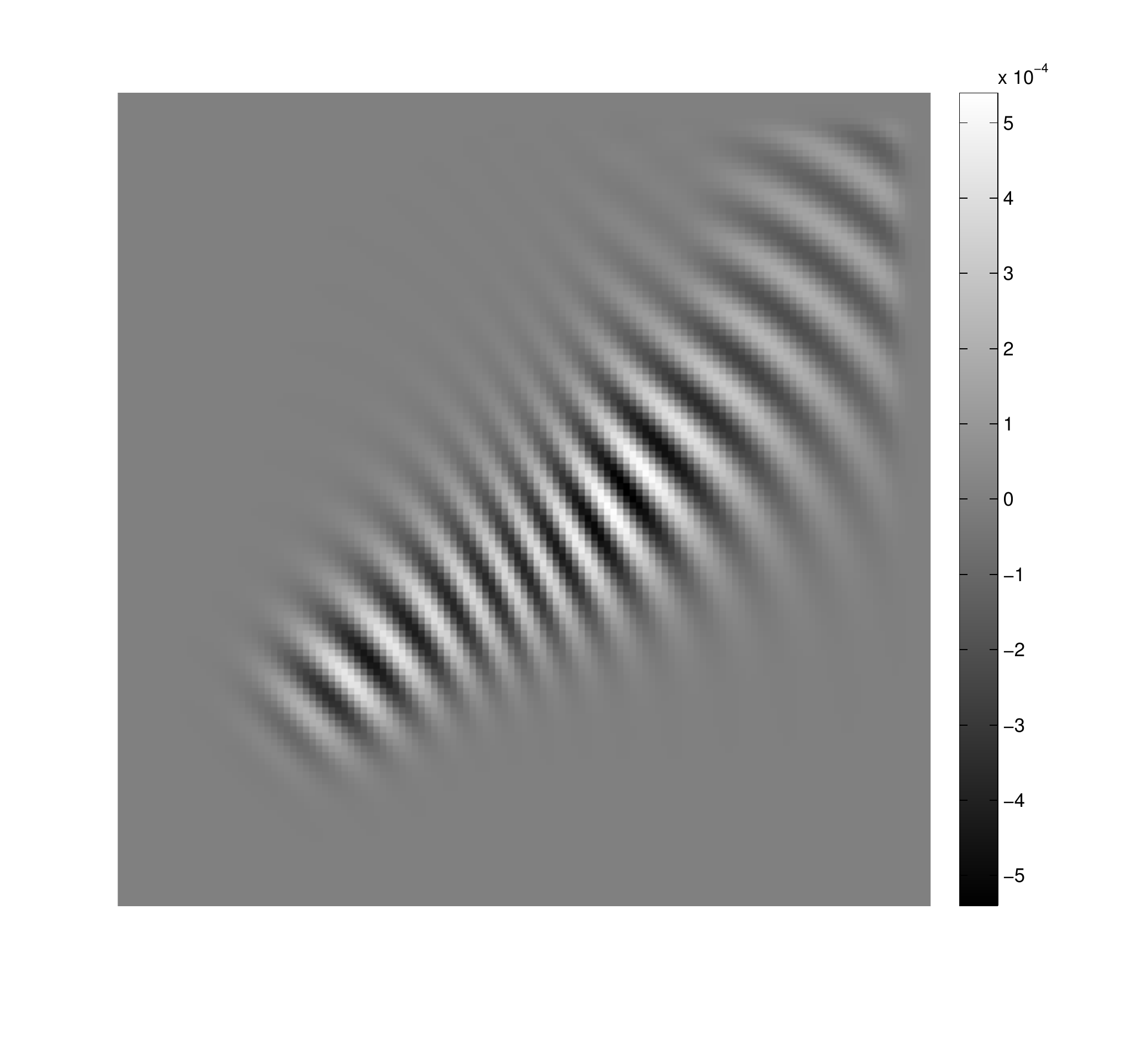}
      \put(42,6){$(b)$}
    \end{overpic}
   \begin{tabular}{ccc|ccc|c|cc|ccc}
     \hline
     &&&\multicolumn{3}{|c|}{$T_{\text{setup}}$}&&\multicolumn{2}{|c}{$N_{\text{iter}}$}&\multicolumn{3}{|c}{$T_{\text{solve}}$}\\
     \hline
     $\omega/(2\pi)$&$q$&$N$& NR & R &ratio&$f(x)$& NR & R & NR & R &ratio\\
     \hline
     \multirow{2}{*}{8}&\multirow{2}{*}{8}&\multirow{2}{*}{$63^3$}&\multirow{2}{*}{48.855}&\multirow{2}{*}{18.490}&\multirow{2}{*}{38\%}&$(a)$&3&4&11.249&15.996&142\%\\
     &&&&&&$(b)$&3&5&11.048&19.581&177\%\\
     \hline
     \multirow{2}{*}{16}&\multirow{2}{*}{8}&\multirow{2}{*}{$127^3$}&\multirow{2}{*}{524.61}&\multirow{2}{*}{163.16}&\multirow{2}{*}{31\%}&$(a)$&4&5&152.32&212.59&140\%\\
     &&&&&&$(b)$&3&5&111.71&213.24&191\%\\
     \hline
     \multirow{2}{*}{32}&\multirow{2}{*}{8}&\multirow{2}{*}{$255^3$}&\multirow{2}{*}{6038.8}&\multirow{2}{*}{1319.0}&\multirow{2}{*}{22\%}&$(a)$&5&6&1676.7&2471.9&147\%\\
     &&&&&&$(b)$&4&5&1345.3&2084.6&155\%\\
     \hline
   \end{tabular}
   \caption{Results for velocity field (b) in Figure
     \ref{fig:c}. Solutions with $\omega/(2\pi)=16$ at $x_1=0.5$ are
     presented. }
\end{table}

\begin{table}[h!]
  \centering
  \begin{overpic}
    [width=0.45\textwidth]{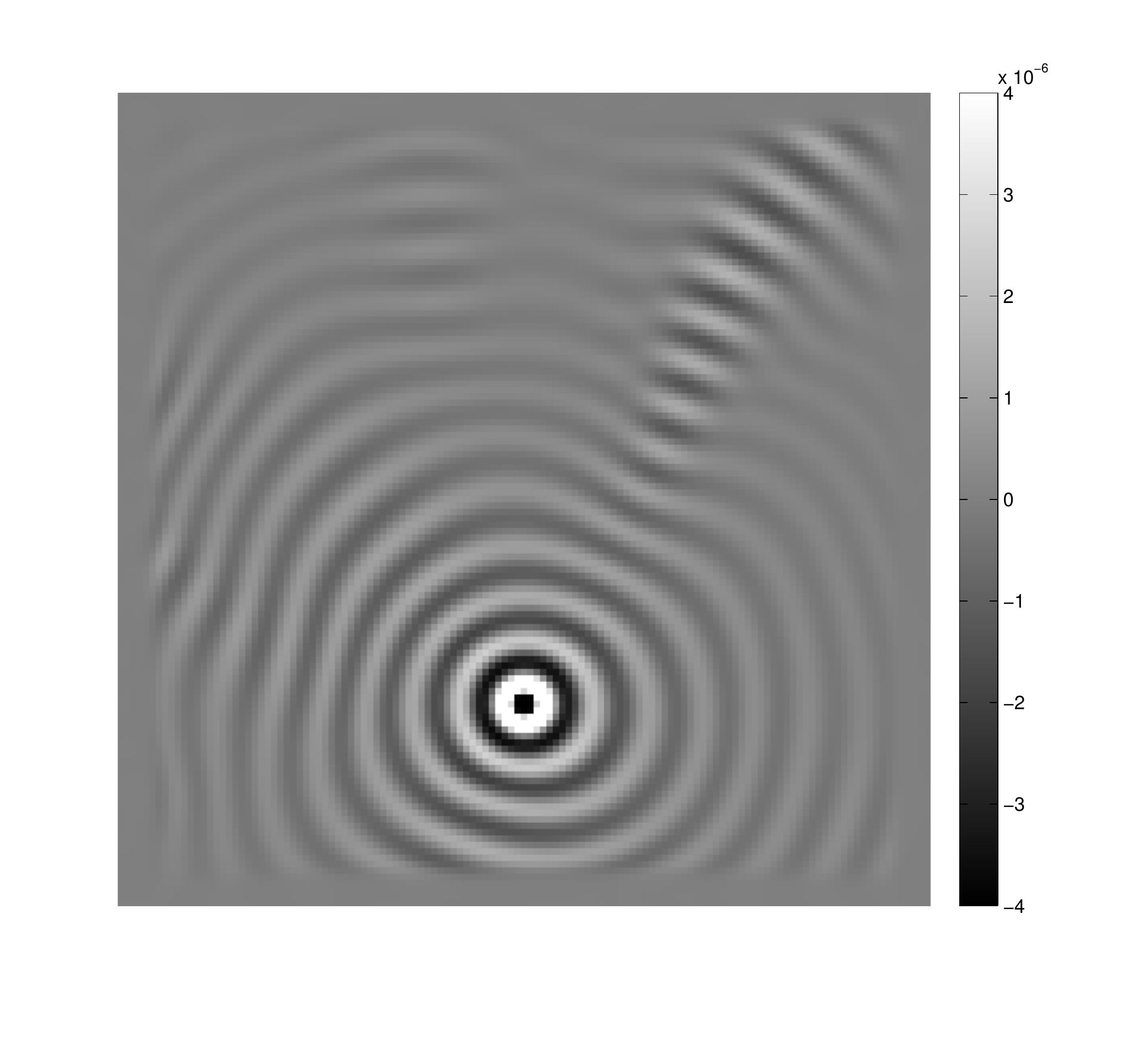}
    \put(42,6){$(a)$}
  \end{overpic}
  \begin{overpic}
    [width=0.45\textwidth]{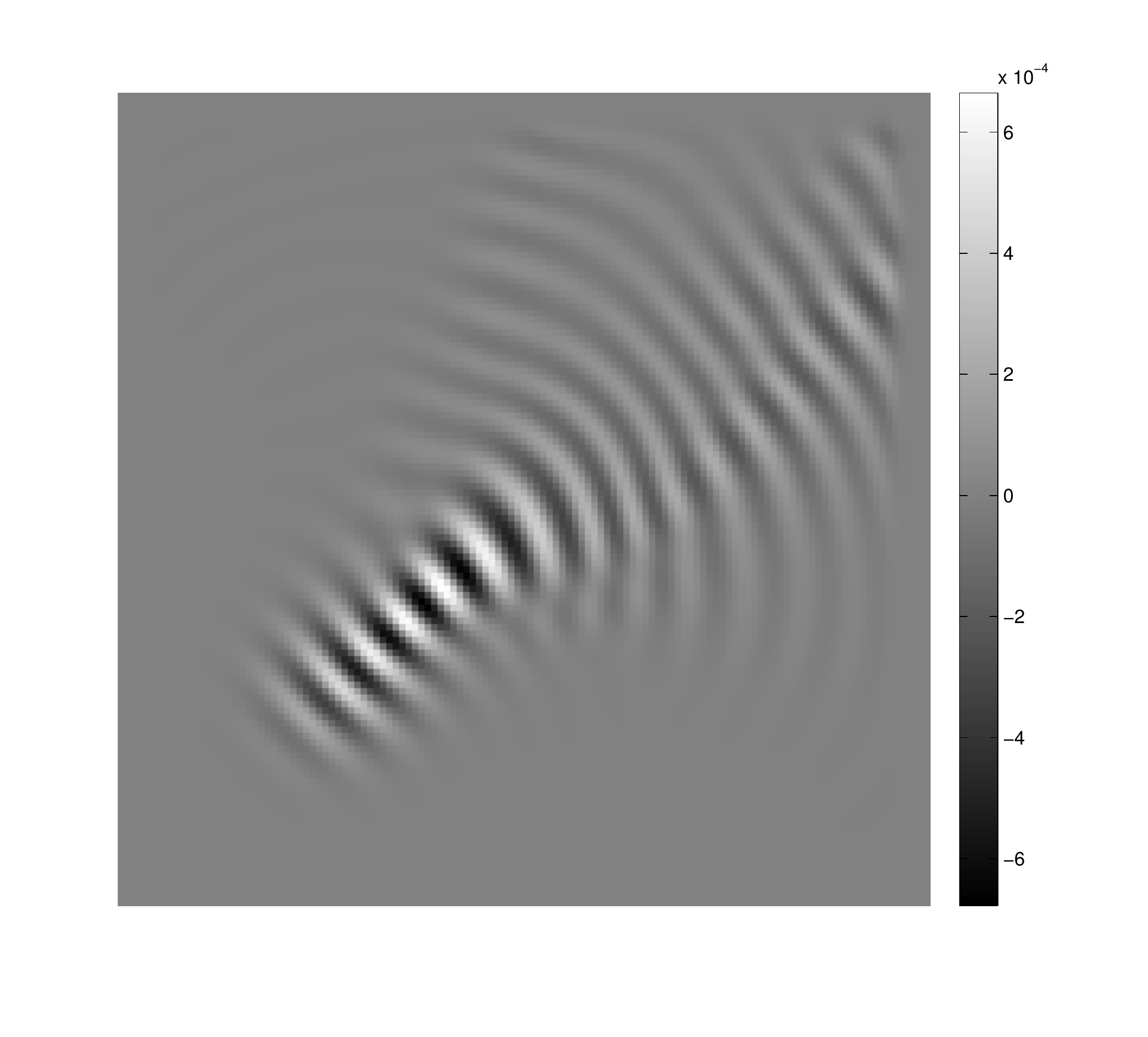}
    \put(42,6){$(b)$}
  \end{overpic}
  \begin{tabular}{ccc|ccc|c|cc|ccc}
    \hline
    &&&\multicolumn{3}{|c|}{$T_{\text{setup}}$}&&\multicolumn{2}{|c}{$N_{\text{iter}}$}&\multicolumn{3}{|c}{$T_{\text{solve}}$}\\
    \hline
    $\omega/(2\pi)$&$q$&$N$& NR & R &ratio&$f(x)$& NR & R & NR & R &ratio\\
    \hline
    \multirow{2}{*}{8}&\multirow{2}{*}{8}&\multirow{2}{*}{$63^3$}&\multirow{2}{*}{48.949}&\multirow{2}{*}{18.213}&\multirow{2}{*}{37\%}&$(a)$&4&4&16.487&16.252&99\%\\
    &&&&&&$(b)$&4&4&15.747&16.513&105\%\\
    \hline
    \multirow{2}{*}{16}&\multirow{2}{*}{8}&\multirow{2}{*}{$127^3$}&\multirow{2}{*}{553.70}&\multirow{2}{*}{147.00}&\multirow{2}{*}{27\%}&$(a)$&4&4&158.41&170.59&108\%\\
    &&&&&&$(b)$&5&5&196.91&215.84&110\%\\
    \hline
    \multirow{2}{*}{32}&\multirow{2}{*}{8}&\multirow{2}{*}{$255^3$}&\multirow{2}{*}{6024.1}&\multirow{2}{*}{1299.2}&\multirow{2}{*}{22\%}&$(a)$&5&5&1599.6&1929.1&121\%\\
    &&&&&&$(b)$&5&6&1635.2&2314.5&142\%\\
    \hline
  \end{tabular}
  \caption{Results for velocity field (c) in Figure
    \ref{fig:c}. Solutions with $\omega/(2\pi)=16$ at $x_1=0.5$ are
    presented. }
\end{table}

From these tests we can make the following observations:
\begin{enumerate}
\item
  The setup time of the recursive preconditioner is significantly dropped
  compared to the non-recursive one. The advantage becomes more and more obvious
  when the problem size gets larger. This is because the setup cost of
  the recursive method scales like $O(N)$ and the non-recursive one scales like
  $O(N^{4/3})$.
\item
  The iteration number of the recursive preconditioner increases only slightly compared to the non-recursive one. Typically it needs
  about 1 more iteration. This is mainly because the recursive method
  solves the quasi-2D auxiliary problems approximately while the
  non-recursive one solves them accurately.
\item
  The application time of the recursive sweeping preconditioner is not
  as fast as the non-recursive method, due to a larger prefactor of
  the time complexity. However, we think this sacrifice is acceptable
  since a huge amount of time is saved in the setup process. One thing
  that needs to be pointed out is that the ratio of the solve time
  increases as the problem size increases, which seems to be
  unexpected since the solve time of the recursive method scales like $O(N)$ and the
  non-recursive one scales like $O(N \log N)$. The reason of this
  behavior is, when the problem size increases, the size of the dense
  linear algebra operations increases as well in the non-recursive
  method since the front size in the nested dissection method gets larger,
  while in the recursive method, the size of the dense linear algebra
  operations in the quasi-1D block LDU setup process and
  solve process is kept the same. Since MATLAB processes large scale
  dense linear algebra operations in a parallel way, the non-recursive
  method gains some advantages from this.
\end{enumerate}

  The memory cost of the recursive method is also advantageous. In our implementation, the recursive method costs only $30\%$ memory compared to the non-recursive method in the $N=255^3$
  case. Theoretically, the memory cost of the recursive method scales
  like $O(N)$ while the non-recursive one scales like
  $O(N\log(N))$. This is another main advantage of the recursive
  method.

\section{Conclusion and Future Work}
\label{sec:conclu}

In this paper, we introduced a new recursive sweeping preconditioner for
the 3D Helmholtz equation based on the moving PML sweeping
preconditioner proposed in \cite{sweeppml}. The idea of the sweeping
preconditioner is used recursively for the auxiliary quasi-2D
problems. Both the setup cost and application cost of the
preconditioner are reduced to strict linear complexity. The iteration
number remains essentially independent of the problem size when
combined with the standard GMRES solver. Numerical results show that
the setup time drops significantly compared to the non-recursive
method, while the solve cost increases only slightly.

Several questions still remain open and some potential improvements
can be made. First, we use the PML to simulate the Sommerfeld
condition. Many other simulations of the absorbing boundary condition
can be implemented and the recursive sweeping idea can be used as long
as the stencil of the simulation is local. Second, the numerical
scheme used in this paper is the standard central difference scheme,
whose dispersion relationship is a poor approximation of the true
one. More accurate numerical schemes can be implemented and the
iteration number may be dropped potentially benefiting from the
increment of the accuracy of the numerical scheme.

Parallel processing can also be introduced to the current recursive
method. First, When sweeping from both sides of the domain, either in
the outer loop of the algorithm or in the inner loop, the processing
of the two fronts can be paralleled so in total it could be $4$ times
faster with parallelization theoretically. Second, the quasi-1D
problems are solved by the block LDU factorization in the current
setting. If instead, we use the 1D nested dissection algorithm for the
quasi-1D problems, then it can be easily paralleled and the total cost
will remain essentially the same. Last, one can notice that, the setup
process of the algorithm is essentially $O(n^2)$ quasi-1D subproblems
which are independent with each other so this process can be done in
parallel, and compared to the original method, which contains only
$O(n)$ quasi-2D independent subproblems, the potential advantages of
parallelization in the setup stage is more obvious here.

There are also several other advantages of the recursive sweeping
method that concern flexibility. First, as mentioned above, the setup
process contains $O(n^2)$ quasi-1D independent subproblems. So if the
velocity field is modified on a subdomain which involves only limited
subproblems, then the factorization can be updated with only a slight
modification on these involved subproblems. Compared to the original
method, where the subproblems are $O(n)$ quasi-2D plates, the
recursive method is more flexible on updating the factorization. This
could be advantageous in seismic imaging where the velocity field is
tested and modified frequently. Second, when the factorization for the
$O(n^2)$ subproblems is done, there are naturally two ways of using
the factorization. One is, as mentioned in this paper, sweeping along
the $x_3$ direction in the outer loop, and sweeping along the $x_2$
direction in the inner loop. Another choice is to do the opposite,
which is sweeping along the $x_2$ direction in the outer loop and the
$x_3$ direction in the inner loop. Each of these two choices shows
some ``bias'' since the residual of the system is accumulated in some
``chosen'' order. So one may ask that, is it possible to combine the
two choices together to make the solve process more flexible such that
the total solve time can be even less? This is another interesting
question to be examined.

\section*{Acknowledgments}
This work was partially supported by the National Science Foundation
under award DMS-1328230 and the U.S. Department of Energy's Advanced
Scientific Computing Research program under award
DE-FC02-13ER26134/DE-SC0009409. We thank Lenya Ryzhik for providing
computing resources and thank Laurent Demanet and Paul Childs for
helpful discussions.

\bibliographystyle{abbrv}
\bibliography{ref}

\begin{thebibliography}{10}

\bibitem{berenger1994pml}
J.-P. Berenger.
\newblock A perfectly matched layer for the absorption of electromagnetic
  waves.
\newblock {\em J. Comput. Phys.}, 114(2):185--200, 1994.

\bibitem{chen2013sourcetrans}
Z.~Chen and X.~Xiang.
\newblock A source transfer domain decomposition method for {H}elmholtz
  equations in unbounded domain.
\newblock {\em SIAM J. Numer. Anal.}, 51(4):2331--2356, 2013.

\bibitem{chen2013sourcetrans2}
Z.~Chen and X.~Xiang.
\newblock A source transfer domain decomposition method for {H}elmholtz
  equations in unbounded domain {P}art {II}: {E}xtensions.
\newblock {\em Numer. Math. Theory Methods Appl.}, 6(3):538--555, 2013.

\bibitem{chew1994pml}
W.~C. Chew and W.~H. Weedon.
\newblock A {3D} perfectly matched medium from modified {M}axwell's equations
  with stretched coordinates.
\newblock {\em Microw. Opt. Techn. Let.}, 7(13):599--604, 1994.

\bibitem{duff1983multifrontal}
I.~S. Duff and J.~K. Reid.
\newblock The multifrontal solution of indefinite sparse symmetric linear
  equations.
\newblock {\em ACM Trans. Math. Software}, 9(3):302--325, 1983.

\bibitem{sweephmf}
B.~Engquist and L.~Ying.
\newblock Sweeping preconditioner for the {H}elmholtz equation: hierarchical
  matrix representation.
\newblock {\em Comm. Pure Appl. Math.}, 64(5):697--735, 2011.

\bibitem{sweeppml}
B.~Engquist and L.~Ying.
\newblock Sweeping preconditioner for the {H}elmholtz equation: moving
  perfectly matched layers.
\newblock {\em Multiscale Model. Simul.}, 9(2):686--710, 2011.

\bibitem{advances}
Y.~A. Erlangga.
\newblock Advances in iterative methods and preconditioners for the {H}elmholtz
  equation.
\newblock {\em Arch. Comput. Methods Eng.}, 15(1):37--66, 2008.

\bibitem{why}
O.~G. Ernst and M.~J. Gander.
\newblock Why it is difficult to solve {H}elmholtz problems with classical
  iterative methods.
\newblock In {\em Numerical analysis of multiscale problems}, volume~83 of {\em
  Lect. Notes Comput. Sci. Eng.}, pages 325--363. Springer, Heidelberg, 2012.

\bibitem{ailu}
M.~J. Gander and F.~Nataf.
\newblock A{ILU} for {H}elmholtz problems: a new preconditioner based on the
  analytic parabolic factorization.
\newblock {\em J. Comput. Acoust.}, 9(4):1499--1506, 2001.

\bibitem{george1973nested}
A.~George.
\newblock Nested dissection of a regular finite element mesh.
\newblock {\em SIAM J. Numer. Anal.}, 10:345--363, 1973.
\newblock Collection of articles dedicated to the memory of George E. Forsythe.

\bibitem{johnson2008pmlnotes}
S.~G. Johnson.
\newblock Notes on perfectly matched layers ({PML}s).
\newblock {\em Lecture notes, Massachusetts Institute of Technology,
  Massachusetts}, 2008.

\bibitem{liu1992multifrontal}
J.~W.~H. Liu.
\newblock The multifrontal method for sparse matrix solution: theory and
  practice.
\newblock {\em SIAM Rev.}, 34(1):82--109, 1992.

\bibitem{parallelsweep}
J.~Poulson, B.~Engquist, S.~Li, and L.~Ying.
\newblock A parallel sweeping preconditioner for heterogeneous 3{D} {H}elmholtz
  equations.
\newblock {\em SIAM J. Sci. Comput.}, 35(3):C194--C212, 2013.

\bibitem{stolk2013domaindecomp}
C.~C. Stolk.
\newblock A rapidly converging domain decomposition method for the {H}elmholtz
  equation.
\newblock {\em J. Comput. Phys.}, 241(0):240 -- 252, 2013.

\bibitem{sweepemfem}
P.~Tsuji, B.~Engquist, and L.~Ying.
\newblock A sweeping preconditioner for time-harmonic {M}axwell's equations
  with finite elements.
\newblock {\em J. Comput. Phys.}, 231(9):3770--3783, 2012.

\bibitem{sweepspectral}
P.~Tsuji, J.~Poulson, B.~Engquist, and L.~Ying.
\newblock Sweeping preconditioners for elastic wave propagation with spectral
  element methods.
\newblock {\em ESAIM Math. Model. Numer. Anal.}, 48(2):433--447, 2014.

\bibitem{sweepem}
P.~Tsuji and L.~Ying.
\newblock A sweeping preconditioner for {Y}ee's finite difference approximation
  of time-harmonic {M}axwell's equations.
\newblock {\em Front. Math. China}, 7(2):347--363, 2012.

\bibitem{vion2014doublesweep}
A.~Vion and C.~Geuzaine.
\newblock Double sweep preconditioner for optimized schwarz methods applied to
  the {H}elmholtz problem.
\newblock {\em J. Comput. Phys.}, 266(0):171 -- 190, 2014.

\bibitem{demanet}
L.~{Zepeda-N{\'u}{\~n}ez} and L.~{Demanet}.
\newblock {The method of polarized traces for the 2D Helmholtz equation}.
\newblock {\em ArXiv e-prints}, Oct. 2014.

\end{thebibliography}
\end{document}